\documentclass[11pt,reqno]{amsart}
\usepackage{amssymb}
\theoremstyle{plain}
\newtheorem{theorem}{Theorem}[section]
\theoremstyle{remark}
\newtheorem{remark}[theorem]{Remark}
\newtheorem{example}[theorem]{Example}
\theoremstyle{plain}
\newtheorem{corollary}[theorem]{Corollary}
\newtheorem{lemma}[theorem]{Lemma}
\newtheorem{proposition}[theorem]{Proposition}

\newtheorem{hypothesis}[theorem]{Hypothesis}
\numberwithin{equation}{section}

\begin{document}
\title[SYMMETRIC OU GENERATORS]
{{\bf SYMMETRIC ORNSTEIN-UHLENBECK SEMIGROUPS}\\
{\bf AND THEIR GENERATORS}}
\author{Anna Chojnowska-Michalik}
\address{Institute of Mathematics, \L\'od\'z University, 
Stefana Banacha 22, 90238 \L\'od\'z, Poland}
\email{katprob@krysia.uni.lodz.pl}
\author{Beniamin Goldys}
\address{School of Mathematics, The University of New South 
Wales, Sydney 2052, Australia}
\email{B.Goldys@unsw.edu.au}
%\today
\thanks{This work was partially supported by the Small 
ARC Grant Scheme}
\begin{abstract}
We provide necessary and sufficient conditions for a  
Hilbert 
space-valued Ornstein-Uhlenbeck process to be reversible with 
respect to its invariant measure $\mu$. For a reversible process 
the domain of its generator in $L^p(\mu )$ is characterized in terms 
of appropriate Sobolev spaces thus extending the Meyer 
equivalence of norms to any symmetric 
Ornstein-Uhlenbeck operator. We provide also a formula for 
the size of the spectral gap of the generator. Those 
results are applied to study the 
Ornstein-Uhlenbeck process in a chaotic environment. 
Necessary and sufficient conditions for a transition semigroup 
$\left(R_t\right)$ to be compact, Hilbert-Schmidt and strong Feller are  
given in terms of the coefficients of the 
Ornstein-Uhlenbeck operator. We show also that the 
existence of spectral gap implies a smoothing property 
of $R_t$  and provide an estimate for the (appropriately 
defined) gradient of $R_t\phi$. Finally, 
in the Hilbert-Schmidt case, we show that for any $\phi\in L^p(\mu 
)$ the 
function $R_t\phi$ is an (almost) classical solution of a 
version of the  
Kolmogorov 
equation. 
\end{abstract}
\maketitle
\tableofcontents
\section{INTRODUCTION}
Consider a linear stochastic differential 
equation 
\begin{equation}\left\{\begin{array}{l}
dZ=AZdt+\sqrt {Q}dW,\\
Z(0)=x,\end{array}
\right.\label{01}\end{equation}
in a separable real Hilbert space $H$. We assume that $Q$ is a 
bounded selfadjoint and nonnegative operator on $H$ and $A$ 
generates on $H$ a strongly continuous semigroup $\left(S(t)\right
)$. 
The process $W$ is a standard cylindrical Wiener process 
on $H$. Under appropriate 
assumptions (see Hypothesis \ref{H}\ below) the solution 
to (\ref{01}), is given by the formula 
\[Z(t,x)=S(t)x+\int_0^tS(t-s)\sqrt {Q}dW(s),\quad t\ge 0.\]
The process $Z$, called the Ornstein-Uhlenbeck process, is Gaussian 
and Markov in $H$ with the transition semigroup 
\[R_t\phi (x)=E\phi\left(Z(t,x)\right),\]
where $\phi$ is a bounded Borel function on $H$. An important 
class of the Ornstein-Uhlenbeck processes are the so 
called reversible Ornstein-Uhlenbeck processes which 
arise in the theory of Interacting Particle Systems and 
other areas of Mathematical Physics. Let us recall that 
a probability measure $\mu$ is said to be symmetrizing for the 
semigroup $\left(R_t\right)$, or equivalently the process $Z$ is 
said to be $\mu$-reversible, if 
\begin{equation}\int_H\psi (x)R_t\phi (x)\mu (dx)=\int_H\phi (x)R_
t\psi (x)\mu (dx),\label{rev}\end{equation}
for bounded Borel functions $\phi ,\psi$. 
If such a measure $\mu$ exists then it is necessarily 
invariant for the semigroup $\left(R_t\right)$. In that case $\left
(R_t\right)$ 
extends to a strongly continuous 
semigroup of contractions on $L^p(H,\mu )$ (still denoted by 
$\left(R_t\right)$) for all $p\in [1,\infty )$. Moreover, (\ref{rev}) implies that 
$R_t$ is symmetric in $L^2(H,\mu )$ for each $t\ge 0$. 
\par
The aim of this paper is to provide necessary and 
sufficient conditions for the invariant measure of an 
arbitrary Ornstein-Uhlenbeck process $Z$ to be 
symmetrizing and to study some important properties 
of a symmetric semigroup $\left(R_t\right)$ and its generator $L$ in 
the spaces $L^p(H,\mu )$ under the sole assumption of existence of 
a nondegenerate invariant measure $\mu$. 
\par
The main idea of this paper may be described as follows. 
Let $H_Q=Q^{1/2}(H)$ with the norm $|x|_Q=\left|Q^{-1/2}x\right|$. We will 
show that  
the semigroup $\left(R_t\right)$ is symmetric in $L^2(H,\mu )$ if and only if 
$H_Q$ is invariant for $\left(S(t)\right)$ and its restriction to $
H_Q$ 
defines a $C_0$-semigroup of symmetric contractions. This fact 
allows us to provide explicit criteria in 
terms of $A$ and $Q$ for the various interesting properties 
of the Ornstein-Uhlenbeck semigroup $\left(R_t\right)$ and its 
generator $L$. In particular, we characterize the domain of $L$ in 
$L^p(H,\mu )$, the spectral gap 
property, compactness, Hilbert-Schmidt property and the 
strong Feller property, and finally, the existence of 
''almost classical''  
solutions to the associated Kolmogorov equation.  
Let us note that the existing 
characterizations  
of those properties for a general Ornstein-Uhlenbeck 
semigroup are usually not easily applicable, see \cite{dz1}, 
\cite{fock}, \cite{reg}, \cite{symp}.  
\par
Let us emphasise that our starting point is the stochastic 
differential equation (\ref{01}) and the associated 
transition semigroup $\left(R_t\right)$. For another approach, where 
the starting point is the Gaussian measure 
$\mu =N\left(0,Q_{\infty}\right)$ and the associated Dirichlet form 
\[\mathcal E(\phi ,\psi )=\frac 12\int_H\left\langle Q^{1/2}D\phi 
(x),Q^{1/2}D\psi (x)\right\rangle\mu (dx),\]
see for example \cite{bogachev}. 
\par
We will describe now the results of this paper in more 
detail. 
\par
In Section 2 we show that the process $Z$ is reversible if 
and only if for every $x\in\mbox{\rm dom}\left(A^{*}\right)$
\begin{equation}Qx\in\mbox{\rm dom}(A)\quad\mbox{\rm and}\quad AQ
x=QA^{*}x.\label{00}\end{equation}
This result was proved in \cite{zabczyk} for $Q=I$. 
Furthermore we show that the operators $S_Q(t)=Q^{-1/2}S(t)Q^{1/2}$ 
are bounded in $H$ and define a $C_0$-semigroup of symmetric 
contractions on $H$. As a consequence, we find that  
$R_t=\mathcal U\Gamma\left(S_Q(t)\right)\mathcal U^{*}$, where $\Gamma\left
(S_Q(t)\right)$ stands for 
the second quantization of the operator $S_Q(t)$ and 
$\mathcal U:L^2(H,\mu )\to L^2(H,\mu )$ is an isometric isomorphism. 
These 
characterizations allow us to express various properties of 
$\left(R_t\right)$ in terms of analogous properties of $\left(S_Q
(t)\right)$. 
\par
In Section 3 we study the domain of $L$ in $L^p(H,\mu )$. The generator 
$L$ of the semigroup $\left(R_t\right)$ 
may be easily evaluated on a dense set 
of cylindrical functions, see \cite{fock}:  
\[L\phi (x)=\frac 12\mbox{\rm tr}\left(QD^2\phi (x)\right)+\left\langle 
x,A^{*}D\phi (x)\right\rangle ,\quad x\in H,\]
where $D$ stands for the Fr\'echet derivative of the 
function $\phi :H\to\mathbb R$. The problem of an explicit 
characterization of the domain of $L$ in $L^p(H,\mu )$ 
was an object of 
intense study for some time, see \cite{shi}, \cite{regd}, \cite{dap}, 
\cite{dapgo}. We use the results from our recent work 
\cite{paley} to give a complete characterization of the 
domain $\mbox{\rm dom}_p(L)$, $p\in (1,\infty )$, of the selfadjoint (in $
L^2(H,\mu )$) 
generator $L$ acting in $L^p(H,\mu )$ in terms 
of appropriately defined Gauss-Sobolev spaces. 
\par
In Section 4 we show that for a symmetric semigroup 
$\left(R_t\right)$, the spectral gap of $L$, and the compactness of $\left
(R_t\right)$ 
are 
determined by the corresponding properties of the 
semigroup $\left(S_Q(t)\right)$ in $H$. In particular, the spectral gap 
of $L$ in $L^2(H,\mu )$ is the same as the spectral gap of $A_Q$, 
the generator 
of $\left(S_Q(t)\right)$ in $H$. Next, we provide necessary 
and sufficient conditions for the strong Feller 
property of the semigroup $\left(R_t\right)$. Finally, we show that 
for a bounded function $\phi$ 
\[\left\|Q^{1/2}DR_t\phi\right\|_{\infty}\le\frac c{\sqrt {t}}\left
\|\phi\right\|_{\infty},\]
for an arbitrary symmetric Ornstein-Uhlenbeck 
semigroup $\left(R_t\right)$ with the spectral gap property. 
\par
In Section 5 we provide necessary and sufficient 
conditions for the semigroup $\left(R_t\right)$ to be of 
Hilbert-Schmidt type. Let us note that this class of 
semigroups includes an important class of strongly 
Feller semigroups, see \cite{reg}. We show that $\left(R_t\right)$ is 
of Hilbert-Schmidt type if and only if the semigroup 
$\left(S_Q(t)\right)$ is exponentially stable and of Hilbert-Schmidt 
type. Subsequently we prove that the latter condition is 
satisfied if $\mu\left(H_Q\right)=1$. As a 
consequence we find that the function $u(t,x)=R_t\phi (x)$ is 
an (almost) classical solution of the Kolmogorov equation
\[\left\{\begin{array}{ll}
\frac {\partial u}{\partial t}(t,x)=\frac 12\mbox{\rm tr}\left(\left
(D^Q\right)^2u(t,x)\right)+\left\langle Q^{-1/2}x,A_QD^Qu(t,x)\right
\rangle ,&t>0,x\in Q^{1/2}(H),\\
u(0,x)=\phi (x),&\phi\in L^p(H,\mu ),\end{array}
\right.\]
where $D^Q$ denotes the Fr\'echet derivative in the 
direction of the subspace $Q^{1/2}(H)$. 
\par
In Section 6 we discuss some examples. In particular, we 
consider an Ornstein-Uhlenbeck process in a chaotic 
environment and provide a detailed analysis of its 
invariant measures and the Spectral Gap Property. 
\par
In the last part of this section we formulate the main 
assumption of the paper. Let 
\begin{equation}Q_t=\int_0^tS(s)QS^{*}(s)ds,\quad t\le\infty .\label{qt}\end{equation}
The following hypothesis is a standing assumption for 
the rest of this paper and the results will be enunciated 
without further recalling it. 
\begin{hypothesis}\label{H}
We assume that 
\[\int_0^{\infty}\mbox{\rm tr}\left(S(s)QS^{*}(s)\right)ds<\infty 
,\]
and the operator $Q_{\infty}$ is injective. 
\end{hypothesis}
If 
Hypothesis \ref{H} holds then the solution to (\ref{01}) is 
well defined and there exists an  
invariant measure $\mu$ for the process $Z$ which is a 
centered Gaussian measure with the covariance operator 
$Q_{\infty}$, see \cite{dz1}. 
\section{CHARACTERIZATION OF SYMMETRIC OU 
SEMIGROUPS}
We will study equation (\ref{01}) in a separable real 
Hilbert space $H$ with the norm $|\cdot |$. 
\subsection{General OU Process}
The next two 
lemmas summarise some basic properties of an arbitrary 
Ornstein-Uhlenbeck semigroup proved in \cite{fock} and 
\cite{reg} 
which will be useful in the sequel. Lemma \ref{a0}\ seems 
to be new. 
\begin{lemma}\label{fock}
(a) We have $S(t)Q_{\infty}^{1/2}(H)\subset Q_{\infty}^{1/2}(H)$ for each $
t\ge 0$. The family 
of operators $S_0(t)=Q_{\infty}^{-1/2}S(t)Q_{\infty}^{1/2}$, $t\ge 
0$, defines a 
$C_0$-semigroup $\left(S_0(t)\right)$ of contractions on $H$. For each 
$t\ge 0$ 
the adjoint $S_0^{*}(t)$  
may be identified with the operator $\overline {Q_{\infty}^{1/2}S^{
*}(t)Q_{\infty}^{-1/2}}$. 
\par\noindent
(b) Moreover, denoting by $A_0$ the generator of the semigroup 
$\left(S_0(t)\right)$, we find that  $K=Q_{\infty}^{1/2}\left(\mbox{\rm dom}\left
(A^{*}\right)\right)$ is a core for for 
$A_0^{*}$  
and 
\begin{equation}\left\langle A_0^{*}h,h\right\rangle =-\frac 12|V
h|^2,\quad\mbox{\rm for}\quad h\in K,\label{21a}\end{equation}
where $V=Q^{1/2}Q_{\infty}^{-1/2}$ with the domain $Q_{\infty}^{1
/2}(H)$. 
\end{lemma}
\begin{proof}
For (a) see Proposition 1 and (a) of Proposition 2 in 
\cite{reg}. 
\par\noindent
(b) If $x\in\mbox{\rm dom}\left(A^{*}\right)$ then $S^{*}(t)x\in\mbox{\rm dom}\left
(A^{*}\right)$ and 
$S_0^{*}(t)Q_{\infty}^{1/2}x=Q_{\infty}^{1/2}S^{*}(t)x$. Hence, $
K\subset\mbox{\rm dom}\left(A_0^{*}\right)$ and 
$S_0^{*}(t)K\subset K$. Since $K$ is dense in $H$, we find that $
K$ is a 
core for $A_0^{*}$ by Theorem 1.9 in \cite{davies}. Putting 
$x=y\in\mbox{\rm dom}\left(A^{*}\right)$ and $h=Q_{\infty}^{1/2}x$ in the Liapunov Equation 
(\ref{61}) below we obtain 
\[2\left\langle A^{*}Q_{\infty}^{-1/2}h,Q_{\infty}^{1/2}h\right\rangle 
=-\left\langle QQ_{\infty}^{-1/2}h,Q_{\infty}^{-1/2}h\right\rangle 
,\]
which yields (\ref{21a}). 
\end{proof}
Let us recall that if $T:H\to H$ is a contraction 
then  
the second quantization $\Gamma (T):L^2(H,\mu )\to L^2(H,\mu )$ of the 
operator $T$ is well defined, see \cite{simon} or \cite{fock} for details. 
\begin{lemma}\label{fock1}
For each $t\ge 0$ we have $R_t=\Gamma\left(S_0^{*}(t)\right)$ and therefore the 
semigroup $\left(R_t\right)$ is symmetric in $L^2(H,\mu )$ if and only if 
the semigroup $\left(S_0(t)\right)$ is symmetric in $H$. If $A_0=
A_0^{*}$ 
then $V$ is closable and $\mbox{\rm dom}\left(\bar V\right)=\mbox{\rm dom}\left
(\sqrt {-A_0}\right)$. 
\end{lemma}
\begin{proof}
Theorem 1 in \cite{fock} yields immediately 
$R_t=\Gamma\left(S_0^{*}(t)\right)$. The second statement follows from Lemma 
2c) in \cite{fock}. To prove the last one, note that 
$A_0=A_0^{*}\le 0$  and thereby the operator $\sqrt {-A_0}$ is well 
defined and closed. Moreover, 
\begin{equation}\left|Vh\right|^2=2\left|\sqrt {-A_0}h\right|^2,\quad 
h\in K,\label{21b}\end{equation}
by (\ref{21a}), with $K$ being a core for $A_0$. Hence, the 
restricted operator $V|K$ is closable and $\mbox{\rm dom}\left(\sqrt {
-A_0}\right)$ is the 
domain of its closure $\overline {V|K}$. Since $K$ is dense in 
$Q_{\infty}^{1/2}(H)=\mbox{\rm dom}(V)$ in the range norm, we get $\overline 
V=\overline {V|K}$ and 
the lemma follows. 
\end{proof}
Let us note the following corollary of Hypothesis \ref{H}. 
\begin{lemma}\label{a0}
We have $\mbox{\rm ker}\left(A_0\right)=\{0\}$ and $\mbox{\rm ker}\left
(A^{*}\right)=\{0\}$. 
\end{lemma}
\begin{proof}
To prove that $\mbox{\rm ker}\left(A_0\right)=\{0\}$ we will show first that 
\begin{equation}\lim_{t\to\infty}S(t)Q_{\infty}^{1/2}x=0,\quad x\in 
H.\label{zero}\end{equation}
Indeed, since by Hypothesis \ref{H} the operators 
$S(s)QS^{*}(s)$ and $S(s)Q_{\infty}S^{*}(s)$ are nonnegative and of trace 
class, we obtain for any orthonormal base $\left\{e_i:i\ge 1\right
\}$ in $H$: 
\[\int_0^t\mbox{\rm tr}\left(S(s)QS^{*}(s)\right)ds=\sum_{i=1}^{\infty}\left
\langle\int_0^tS(s)QS^{*}(s)e_ids,e_i\right\rangle\]
\[=\sum_{i=1}^{\infty}\left\langle\left(Q_{\infty}-S(t)Q_{\infty}
S^{*}(t)\right)e_i,e_i\right\rangle\]
\[=\mbox{\rm tr}\left(Q_{\infty}\right)-\left\|Q_{\infty}^{1/2}S^{
*}(t)\right\|^2_{HS}.\]
Hence, taking into account that 
$\left\|Q_{\infty}^{1/2}S^{*}(t)\right\|_{HS}=\left\|S(t)Q_{\infty}^{
1/2}\right\|_{HS}$ and in view of Hypothesis 
\ref{H} we find that 
\[\mbox{\rm tr}\left(Q_{\infty}\right)=\lim_{t\to\infty}\int_0^t\mbox{\rm tr}\left
(S(s)QS^{*}(s)\right)ds\]
\[=\mbox{\rm tr}\left(Q_{\infty}\right)-\lim_{t\to\infty}\left\|S
(t)Q_{\infty}^{1/2}\right\|_{HS},\]
which in particular implies (\ref{zero}). Finally, if  
$x\in\mbox{\rm ker}\left(A_0\right)$ then $S_0(t)x=x$, hence $S(t
)Q_{\infty}^{1/2}x=Q_{\infty}^{1/2}x$ 
for all $t\ge 0$. Therefore, by (\ref{zero}) $x\in\mbox{\rm ker}\left
(Q_{\infty}\right)=\{0\}$. 
\par\noindent
Similarly, if $x\in\mbox{\rm ker}\left(A^{*}\right)$ then, taking into account that 
\[Q_tx=Q_{\infty}x-S(t)Q_{\infty}S^{*}(t)x,\quad x\in H,\]
we find that 
\[\lim_{t\to\infty}Q_{\infty}^{1/2}S^{*}(t)x=0,\]
which yields $x=0$. 
\end{proof}
\subsection{Characterizations of the Symmetry}
The next theorem provides necessary and sufficient conditions for the 
semigroup $\left(R_t\right)$ to be symmetric in $L^2(H,\mu )$. This problem 
has been solved in \cite{zabczyk} for the case $Q=I$. 
\begin{theorem}\label{t61} The following 
conditions are equivalent. 
\par\noindent
(i) The semigroup $\left(R_t\right)$ is symmetric in $L^2(H,\mu )$. 
\par\noindent
(ii) If $x\in\mbox{\rm dom}\left(A^{*}\right)$ then $Qx\in\mbox{\rm dom}
(A)$ and 
\begin{equation}AQx=QA^{*}x.\label{60}\end{equation}
(iii) $S(t)Q=QS^{*}(t)$ for all $t\ge 0$. 
\end{theorem}
\begin{proof}
$(i)\implies (ii)$. By Lemma \ref{fock1}\ $R_t$ is 
symmetric if and only if $A_0^{*}=A_0$. Therefore 
$S_0^{*}(t)=S_0(t)$ for all $t\ge 0$ or, equivalently, 
\begin{equation}S(t)Q_{\infty}=Q_{\infty}S^{*}(t).\label{60b}\end{equation}
If Hypothesis \ref{H} holds then 
\begin{equation}\left\langle A^{*}x,Q_{\infty}y\right\rangle +\left
\langle A^{*}y,Q_{\infty}x\right\rangle =-\langle Qx,y\rangle\label{61}\end{equation}
for $x,y\in\mbox{\rm dom}\,\left(A^{*}\right)$ (see chapter 11.2 of \cite{dz1}). It follows 
from (\ref{61}) that if $x\in\mbox{\rm dom}\,\left(A^{*}\right)$ 
then $Q_{\infty}x\in\mbox{\rm dom}\,(A)$  
and thereby, (\ref{60b}) yields 
\begin{equation}AQ_{\infty}x=Q_{\infty}A^{*}x=-\frac 12Qx,\label{62}\end{equation}
for all $x\in\mbox{\rm dom}\,\left(A^{*}\right)$. 
Take 
$y\in\mbox{\rm dom}\left(\left(A^{*}\right)^2\right)$. Since $x=A^{
*}y\in\mbox{\rm dom}\left(A^{*}\right)$ we conclude 
from (\ref{62}) that 
\begin{equation}Q_{\infty}A^{*}y\in\mbox{\rm dom}(A),\label{62a}\end{equation}
and 
\begin{equation}AQ_{\infty}A^{*}y=-\frac 12QA^{*}y.\label{62b}\end{equation}
By (\ref{62}) we also have 
\begin{equation}Q_{\infty}A^{*}y=-\frac 12Qy,\label{62c}\end{equation}
which combined with (\ref{62a}) implies that $Qy\in\mbox{\rm dom}
(A)$ 
and 
\begin{equation}AQ_{\infty}A^{*}y=-\frac 12AQy.\label{62d}\end{equation}
From (\ref{62b}) and (\ref{62d}) we obtain 
\begin{equation}AQx=QA^{*}x,\quad x\in\mbox{\rm dom}\left(\left(A^{
*}\right)^2\right).\label{62e}\end{equation}
Since $\mbox{\rm dom}\left(\left(A^{*}\right)^2\right)$ is a core for $
A^{*}$ and the right hand 
side of (\ref{62e}) is well defined for $x\in\mbox{\rm dom}\left(
A^{*}\right)$, (ii) 
follows. 
\par\noindent
{\em (ii)$\implies$(iii)\/} By assumption $\left(\lambda -A\right
)Qy=Q\left(\lambda -A^{*}\right)y$ for 
$\lambda\in\mathbb R$ and $y\in\mbox{\rm dom}\left(A^{*}\right)$. Hence, for a certain $
\lambda_0$ 
and all $\lambda >\lambda_0$ 
\[\left(\lambda -A\right)^{-1}Qx=Q\left(\lambda -A^{*}\right)^{-1}
x,\quad x\in H.\]
Then using the formula for the resolvent of generator 
and properties of the Laplace transform we obtain 
(iii). 
\par\noindent
{\em (iii)$\implies$(i)\/} From Hypothesis \ref{H} and (iii) 
\[S(t)Q_{\infty}x=\int_0^{\infty}S(t+s)QS^{*}(s)xds\]
\[=\int_0^{\infty}S(s)QS^{*}(t+s)xds=Q_{\infty}S^{*}(t)x,\]
for $x\in H$, which yields $S_0(t)=S_0^{*}(t)$ and (i) follows 
from Lemma \ref{fock1}.
\end{proof}
If the Ornstein-Uhlenbeck process (\ref{01}) is diagonal, 
that is there exists a joint eigenbasis  
\[Ae_k=\alpha_ke_k,\quad Qe_k=q_ke_k,\quad k\ge 1,\]
for $A$ and $Q$, then obviously the corresponding Ornstein-Uhlenbeck 
semigroup is symmetric provided Hypothesis \ref{H} 
holds. It is easy 
to see that an Ornstein-Uhlenbeck process with the  
symmetric transition semigroup 
need not be diagonal. As the simplest example 
it is enough to take $H=\mathbb R^2$, 
\[A=\left(\begin{array}{cc}
a&c\\
d&b\end{array}
\right)\quad\mbox{\rm and}\quad Q=\left(\begin{array}{cc}
1&0\\
0&q\end{array}
\right)\]
with $0<q\neq 1$. Then Hypothesis \ref{H} holds 
and the corresponding Ornstein-Uhlenbeck semigroup is 
symmetric if and only if 
\[a<0,\quad\mbox{\rm det}(A)>0,\quad d=cq,\quad\mbox{\rm and}\quad 
(a-b)^2+4c^2q>0.\]
\begin{corollary}\label{new1}
Let $\left(R_t\right)$ be symmetric. Then the following holds. 
\par\noindent
(i) $\mbox{\rm im}\left(Q_{\infty}\right)\subset\mbox{\rm dom}(A)$ and the operator 
$AQ_{\infty}=-\frac 12Q$ is bounded, symmetric and negative. 
\par\noindent
(ii) $\mbox{\rm im}\left(Q\right)\subset\mbox{\rm im}(A)$. 
\par\noindent
(iii) If  
\begin{equation}\mbox{\rm ker}\left(A\right)=\{0\},\label{29a}\end{equation}
then 
\[Q_{\infty}=-\frac 12A^{-1}Q=-\frac 12\overline {Q\left(A^{*}\right
)^{-1}}.\]
\end{corollary}
\begin{proof}\ 
The proof follows immediately from (\ref{62}). 
\end{proof}
\begin{remark}
Note that (\ref{29a}) holds if the semigroup $\left(S(t)\right)$ is 
stable. In particular, if $H=\mathbb R^d$, then (\ref{29a}) is 
a consequence of Hypothesis \ref{H}. Indeed, by Theorem 
\ref{t61a}(i) below the operator $Q^{-1}$ is bounded in the case of 
$H=\mathbb R^d$ hence, 
\[\int_0^{\infty}\left\|S(t)\right\|^2dt<\infty\]
by Hypothesis \ref{H}, which implies the exponential 
stability of the semigroup $\left(S(t)\right)$. 
In Section 6 we 
provide an example in which $\left(R_t\right)$ is symmetric and has 
the Spectral Gap Property but (\ref{29a}) is not satisfied. 
\end{remark}
\begin{theorem}\label{t61a}
The semigroup $\left(R_t\right)$ is symmetric in $L^2(H,\mu )$ if and only 
if the following conditions are satisfied. 
\par\noindent
(i) $Q$ is injective. 
\par\noindent
(ii) We have 
\begin{equation}S(t)Q^{1/2}(H)\subset Q^{1/2}(H),\quad t\ge 0.\label{60a}\end{equation}
(iii) The family of operators 
\[S_Q(t)=Q^{-1/2}S(t)Q^{1/2},\quad t\ge 0,\]
defines a symmetric $C_0$-semigroup of contractions on $H$. 
\end{theorem}
\begin{proof}
Let $\left(R_t\right)$ be symmetric. Suppose that $Qx=0$. Then by  
Hypothesis \ref{H} and (iii) of Theorem \ref{t61} 
\[Q_{\infty}x=\int_0^{\infty}S(2t)Qxdt=0,\]
and since $Q_{\infty}$ is injective (i) follows. It follows from 
(\ref{62b}), (\ref{62d}) and (i) of Corollary \ref{new1}\ 
that 
\begin{equation}-QA^{*}x=2AQ_{\infty}A^{*}x=-AQx,\quad x\in\mbox{\rm dom}\left
(A^{*}\right),\label{65a}\end{equation}
and therefore (\ref{62}) implies 
\begin{equation}\left\langle -QA^{*}x,x\right\rangle\ge 0,\quad x
\in\mbox{\rm dom}\left(A^{*}\right).\label{65}\end{equation}
To prove that $S(t)Q^{1/2}(H)\subset Q^{1/2}(H)$ we will show that 
\begin{equation}\left|Q^{1/2}S^{*}(t)x\right|^2\le\left|Q^{1/2}x\right
|^2,\quad x\in H.\label{66}\end{equation}
To this end note that for $x\in\mbox{\rm dom}\left(A^{*}\right)$ (\ref{62}) implies 
\begin{equation}Qx=-2\int_0^{\infty}S(s)QS^{*}(s)A^{*}xds=2\int_0^{
\infty}S(s)\left(-QA^{*}\right)S^{*}(s)xds,\label{67}\end{equation}
and therefore, 
\[\left\langle S(t)QS^{*}(t)x,x\right\rangle =2\int_0^{\infty}\left
\langle S(t+s)\left(-QA^{*}\right)S^{*}(t+s)x,x\right\rangle ds\]
\[=2\int_t^{\infty}\left\langle -QA^{*}S^{*}(u)x,S^{*}(u)x\right\rangle 
du\le\left\langle Qx,x\right\rangle ,\]
where the last inequality follows from (\ref{65}) and 
(\ref{67}). 
For arbitrary $x\in H$ (\ref{66}) follows by the density of $\mbox{\rm dom}\left
(A^{*}\right)$ in 
$H$. Finally, by Proposition B.1 of 
\cite{dz1}, (\ref{66}) yields $S(t)Q^{1/2}(H)\subset Q^{1/2}(H)$. Hence, 
the operator $S_Q(t)=Q^{-1/2}S(t)Q^{1/2}$ is bounded 
on $H$ for each $t\ge 0$. Obviously, operators $S_Q(t)$ satisfy 
the semigroup property and moreover, 
\begin{equation}S_Q^{*}(t)=\overline {Q^{1/2}S^{*}(t)Q^{-1/2}},\quad 
t\ge 0.\label{67a}\end{equation}
Since $\left(R_t\right)$ is symmetric, Theorem \ref{t61}\ yields 
$S(t)Q=QS^{*}(t)$ and therefore $S_Q(t)x=S_Q^{*}(t)x$ for $x\in Q^{
1/2}(H)$ 
and by 
the density argument, 
\begin{equation}S_Q^{*}(t)x=S_Q(t)x,\quad x\in H.\label{67b}\end{equation}
Putting in (\ref{66}) $x=Q^{-1/2}y$ with 
$y\in Q^{1/2}(H)$ we obtain 
\begin{equation}\left|S_Q^{*}(t)y\right|\le |y|,\quad y\in Q^{1/2}
(H).\label{68}\end{equation}
Since $Q^{1/2}(H)$ is dense in $H$, (\ref{68}) holds for all $y\in 
H$. 
The function $t\to S_Q^{*}(t)x$ is weakly continuous at zero for 
every $x\in Q^{1/2}(H)$, hence for $x\in H$, since $\left\|S_Q^{*}
(t)\right\|\le 1$ for 
all $t$. Therefore, the 
semigroup $\left(S_Q^{*}(t)\right)$ is a strongly continuous contraction 
semigroup in $H$ and so is $\left(S_Q(t)\right)$.
\par
Assume now that conditions (i)-(iii) are satisfied and 
let $x=Q^{1/2}y$. Then by symmetry of $\left(S_Q(t)\right)$ we have 
$Q^{-1/2}S(t)Qy=Q^{1/2}S^{*}(t)y$, or equivalently, 
$S(t)Qy=QS^{*}(t)y$, which proves condition (iii) of Theorem 
\ref{t61}. 
\end{proof}
\subsection{The First Consequences of Symmetry}
\begin{corollary}\label{new}
Let $\left(R_t\right)$ be symmetric. Then the generator $A_Q$ of the 
semigroup $\left(S_Q(t)\right)$ is injective. 
\end{corollary}
\begin{proof}
If $x\in\mbox{\rm ker}\left(A_Q\right)$ then $S_Q(t)x=x$ for each $
t\ge 0$ and hence 
$S(t)Q^{1/2}x=Q^{1/2}x$. By Hypothesis \ref{H} 
we find that 
\[\int_0^{\infty}\left|Q^{1/2}x\right|^2dt=\int_0^{\infty}\left|S
(t)Q^{1/2}x\right|^2dt\]
\[\le |x|^2\int_0^{\infty}\left\|S(t)Q^{1/2}\right\|^2_{HS}dt=|x|^
2\int_0^{\infty}\left\|Q^{1/2}S^{*}(t)\right\|^2_{HS}dt<\infty ,\]
which yields $x\in\mbox{\rm ker}(Q)=\{0\}$. 
\end{proof}
\begin{theorem}\label{tnew1}
Assume that $\left(R_t\right)$ is symmetric. Then 
there exists an isometric isomorphism $\mathcal U:L^2(H,\mu )\to 
L^2(H,\mu )$ 
such that 
\begin{equation}R_t=\mathcal U^{-1}\Gamma\left(S_Q(t)\right)\mathcal 
U,\quad t\ge 0.\label{gamma}\end{equation}
\end{theorem}
\begin{proof}
Note first that $K$ is a core for $A_0=A_0^{*}$. Hence it follows 
from (\ref{21a}) and Lemma \ref{fock1}\ that $A_0=-\frac 12V^{*}\bar {
V}$. 
By Lemma \ref{a0}\ $A_0$ is injective, hence so is $\sqrt {-A_0}$ and 
since (\ref{21b}) extends on $\mbox{\rm dom}\left(\sqrt {-A_0}\right
)=\mbox{\rm dom}\left(\bar V\right)$, the 
operator $\bar {V}$ is also injective, in particular 
$\bar {V}^{-1}=\overline {V^{-1}}$. Moreover, $\overline {\mbox{\rm im}
(V)}=H$, because by Theorem 
\ref{t61a}(i) $\mbox{\rm im}(Q)$ is dense in $H$. Therefore, 
the polar 
decomposition of $\bar {V}$ is given by 
\[\bar {V}=U\sqrt {-2A_0},\]
where $U:H\to H$ is an isometric isomorphism. Since 
$Q_{\infty}^{1/2}(H)$ is invariant for $S_0^{*}(t)$, $t\ge 0$, we have for 
$h\in Q^{1/2}(H)$ 
\[S_Q^{*}(t)h=Q^{1/2}Q_{\infty}^{-1/2}S_0^{*}(t)Q_{\infty}^{1/2}Q^{
-1/2}h\]
\[=VS_0^{*}(t)V^{-1}h=U\left|\bar V\right|S_0^{*}(t)\left|\bar V\right
|^{-1}U^{-1}h\]
\begin{equation}=US_0^{*}(t)U^{-1}h,\label{new3}\end{equation}
where the last equality holds because $\left|\bar V\right|=\sqrt {
-2A_0}$ 
commutes with $S_0(t)=S_0^{*}(t)$. By (i) of Theorem \ref{t61a}\ 
$\mbox{\rm im}\left(Q^{1/2}\right)$ is dense in $H$ and therefore (\ref{new3}) holds 
for $h\in H$. By Lemma \ref{fock1}, 
(\ref{new3}) and the properties of the second 
quantization operator $\Gamma$ (see \cite{simon} or \cite{fock}) we obtain 
\[R_t=\Gamma\left(S_0(t)\right)=\Gamma\left(U^{-1}S_Q(t)U\right)=
\Gamma\left(U^{-1}\right)\Gamma\left(S_Q(t)\right)\Gamma\left(U\right
),\]
and putting $\mathcal U=\Gamma (U)$ we complete the proof of 
the theorem. 
\end{proof}
\begin{proposition}\label{trkhst}
Assume that $\left(R_t\right)$ is symmetric. Then for each $t>0$ 
\begin{equation}Q_t^{1/2}(H)=Q^{1/2}\left(\mbox{\rm dom}\sqrt {-A_
Q}\right).\label{rkhs}\end{equation}
\end{proposition}
\begin{proof}
By Corollary B.7 in \cite{dz1} 
$\mbox{\rm im}\left(Q_t^{1/2}\right)=\mbox{\rm im}(T)$, where $T:
L^2(0,t;H)\to H$ is given by 
\[Tu=\int_0^tS(t-s)Q^{1/2}u(s)ds.\]
By Theorem \ref{t61a} $S(t-s)Q^{1/2}u(s)=Q^{1/2}S_Q(t-s)u(s)$ 
and therefore 
\[Tu=Q^{1/2}\int_0^tS_Q(t-s)u(s)ds.\]
By Remarks B8 and A.18 in \cite{dz1} 
\begin{equation}\left\{\int_0^tS_Q(t-s)u(s)ds;u\in L^2(0,t;H)\right
\}=D_{A_Q}\left(\frac 12,2\right)=\mbox{\rm dom}\sqrt {-A_Q}\label{inter}\end{equation}
for exponentially stable and symmetric semigroup 
$\left(S_Q(t)\right)$. Note that for $\lambda >0$ we have 
$\mbox{\rm dom}\left(\lambda -A_Q\right)=\mbox{\rm dom}\left(-A_Q\right
)$ and 
$\mbox{\rm dom}\sqrt {\lambda -A_Q}=\mbox{\rm dom}\sqrt {-A_Q}$ with the equivalent graph norms. 
Since the interpolation space 
$D_{A_Q}\left(\frac 12,2\right)$ is completely determined by the pair of 
spaces $\left(\mbox{\rm dom}\left(A_Q\right),H\right)$ it follows that 
\[D_{A_Q}\left(\frac 12,2\right)=D_{A_Q-\lambda}\left(\frac 12,2\right
).\]
Finally, 
\[\left\{\int_0^te^{-\lambda (t-s)}S_Q(t-s)u(s)ds;u\in L^2(0,t;H)\right
\}\]
\[=\left\{\int_0^tS_Q(t-s)u(s)ds;u\in L^2(0,t;H)\right\},\]
and therefore (\ref{inter}) still holds 
if $\left(S_Q(t)\right)$ is only bounded and symmetric. 
\end{proof}
\begin{remark}\label{t61b}
Let $H_Q$ denote the space $Q^{1/2}(H)$ endowed with the 
norm $|x|_Q=\left|Q^{-1/2}x\right|$. It follows from Theorem \ref{t61a} that $\left
(R_t\right)$ is 
symmetric if and only if the semigroup $\left(S(t)\right)$ restricted 
to $H_Q$ defines a $C_0$-semigroup $\left(S_r(t)\right)$ of symmetric contractions 
in this space. Its generator $A_r$ is the part of $A$ in $H_Q$. By 
Corollary \ref{new} 
$\mbox{\rm ker}\left(A_r\right)=\{0\}$. Moreover, 
\begin{equation}Q_t^{1/2}(H)=\mbox{\rm dom}\left(\sqrt {-A_r}\right
).\label{r1}\end{equation}
Indeed, by Proposition \ref{trkhst} 
$Q_t^{1/2}(H)=Q^{1/2}A_Q^{-1/2}(H)$, hence any $x\in Q_t^{1/2}(H)$ can be 
written (see for example p. 70 of \cite{pazy}) in the form 
\[x=Q^{1/2}\int_0^{\infty}\frac 1{\sqrt {t}}e^{-\lambda t}S_Q(t)h
dt,\quad h\in H,\]
for a certain $\lambda >0$. 
Therefore, 
\begin{equation}x=\int_0^{\infty}\frac 1{\sqrt {t}}e^{-\lambda t}
S(t)Q^{1/2}hdt,\quad h\in H,\label{hs}\end{equation}
with the integral convergent in the norm $\left|\cdot\right|_Q$. This,  
again by the formula on p. 70 of \cite{pazy}, implies that 
$x\in\mbox{\rm dom}\left(\sqrt {-A_r}\right)$ and (\ref{r1}) follows. 
\end{remark}
\section{IDENTIFICATION OF DOMAINS}
Let $K=Q_{\infty}^{1/2}\left(\mbox{\rm dom}\left(A^{*}\right)\right
)$ and let 
\[\mathcal P(K)=\mbox{\rm lin}\left\{\phi_k^n:n\ge 0,k\in K\right
\},\]
where $\phi_k(x)=\left\langle x,Q_{\infty}^{-1/2}k\right\rangle$. By Lemma \ref{fock}\ $
K$ is a 
core for $A_0^{*}$ and by (15) in \cite{fock} $\mathcal P(K)$ is a 
core for $L$. 
We denote by $\mathcal W_{AQ}^{1,p}$ the completion of $\mathcal 
P(K)$ in the 
norm 
\begin{equation}\left\|\phi\right\|_{1,p,AQ}=\left(\left\|\phi\right
\|_p^p+\left\|\left(-AQ\right)^{1/2}D\phi\right\|_p^p\right)^{1/p}
.\label{n1}\end{equation}
By $\mathcal W_Q^{1,p}$ we denote the completion of $\mathcal P(K
)$ with 
respect to the norm $\left\|\cdot\right\|_{1,p,Q}$, where 
\begin{equation}\left\|\phi\right\|_{1,p,Q}^p=\left\|\phi\right\|_
p^p+\int_H\left|Q^{1/2}D\phi (x)\right|^p\mu (dx).\label{n2}\end{equation}
The completion of $\mathcal P(K)$ with respect to the norm 
\begin{equation}\left\|\phi\right\|_{2,p,Q}=\left(\left\|\phi\right
\|_{1,p,Q}^p+\int_H\left\|Q^{1/2}D^2\phi (x)Q^{1/2}\right\|^p_{HS}
\mu (dx)\right)^{1/p}\label{n3}\end{equation}
where $\left\|\cdot\right\|_{HS}$ denotes the Hilbert-Schmidt norm of an 
operator, will be denoted by $\mathcal W_Q^{2,p}$. 
\begin{theorem}\label{t67} Assume that $L$ is selfadjoint. 
Then for every $p\in (1,\infty )$ the spaces $\mathcal W_Q^{1,p}$, 
$\mathcal W_{AQ}^{1,p}$ and $\mathcal W_Q^{2,p}$ may be identified as 
subspaces of $L^p(H,\mu )$. Moreover, 
\begin{equation}\mbox{\rm dom}_p\left(-L\right)^{1/2}=\mathcal W^{
1,p}_Q,\label{s1}\end{equation}
and 
\begin{equation}\mbox{\rm dom}_p\left(L\right)=\mathcal W_Q^{2,p}
\cap\mathcal W_{AQ}^{1,p}.\label{s2}\end{equation}
\end{theorem}
In order to prove this theorem we will recall some facts from 
\cite{paley}. 
We define first the Malliavin 
gradient 
\[D_I=Q_{\infty}^{1/2}D,\quad\mbox{\rm dom}\left(D_I\right)=\mathcal 
P(K).\]
Taking into account that $-A_0$ is 
nonnegative and selfadjoint we define the gradients 
\[D_{A_0}^n=\left(\left(-A_0\right)^{1/2}\right)^{\otimes n}D_I^n
,\quad\mbox{\rm dom}\left(D_{A_0}^n\right)=\mathcal P\left(K\right
),\quad n=1,2.\]
It was shown in \cite{paley} that for $p\in (1,\infty )$ $D_{A_0}$ is 
closable in $L^p(H,\mu )$ and $D_{A_0}^2$ is closable in 
$\mbox{\rm dom}_p\left(\overline {D_{A_0}}\right)$ endowed with the graph norm in $
L^p(H,\mu )$. 
The closed extensions are again denoted by $D_{A_0}^n$, $n=1,2$. 
The next theorem is a special case of Theorem 5.3 in 
\cite{paley}. 
\begin{theorem}\label{td2}
For each $p\in (1,\infty )$, there exist 
$a_p,b_p>0$ such that for $\phi\in\mbox{\rm dom}_p\left(\sqrt {-L}\right
)$ 
\begin{equation}a_p\left(\|\phi\|_p+\left\|D_{A_0}\phi\right\|_p\right
)\le\left\|\sqrt {I-L}\phi\right\|_p\le b_p\left(\|\phi\|_p+\left
\|D_{A_0}\phi\right\|_p\right),\label{m1}\end{equation}
and for $\phi\in\mbox{\rm dom}\left(L\right)$ 
\[a_p\left(\left\|\phi\right\|_p+\left\|\sqrt {I-A_0}D_{A_0}\phi\right
\|_p+\left\|D_{A_0}^2\phi\right\|_p\right)\le\left\|\left(I-L^{}\right
)\phi\right\|_p\]
\begin{equation}\le b_p\left(\left\|\phi\right\|_p+\left\|\sqrt {
I-A_0}D_{A_0}\phi\right\|_p+\left\|D_{A_0}^2\phi\right\|_p\right)
.\label{m2}\end{equation}
\end{theorem}
We will prove now Theorem \ref{t67}
\begin{proof}\ Note first, that by (\ref{65}) 
the operator 
 $AQ=QA^{*}$ with $\mbox{\rm dom}\left(AQ\right)=\mbox{\rm dom}\left
(A^{*}\right)$ may be extended 
in the sense of Friedrichs to a selfadjoint and 
nonpositive operator in $H$.
We will use the same notation $AQ$ for the Friedrichs 
extension of the operator $\left(AQ,\mbox{\rm dom}\left(A^{*}\right
)\right)$. 
\par
In view of Theorem \ref{td2} it remains to identify the 
Sobolev norms given in (\ref{n1}), (\ref{n2}) and (\ref{n3}) 
with the appropriate norms in (\ref{m1}) and (\ref{m2}). 
The 
relationship (\ref{s1}) was proved in \cite{shi}. We 
repeat here the argument using our notation for the sake of 
completeness. 
For $\phi\in\mathcal P(K)$ we have 
\[\left|D_{A_0}\phi (x)\right|^2=\left|\sqrt {-A_0}Q_{\infty}^{1/
2}D\phi (x)\right|^2\]
\[=\left\langle -Q_{\infty}^{1/2}A_0Q_{\infty}^{1/2}D\phi (x),D\phi 
(x)\right\rangle\]
\[=\left\langle -AQ_{\infty}D\phi (x),D\phi (x)\right\rangle =\frac 
12\left|Q^{1/2}D\phi (x)\right|^2,\]
by (i) of Corollary \ref{new1} and (\ref{62}). Hence the operator 
$\left(Q^{1/2}D,\mathcal P(K)\right)$ is closable in $L^p(H,\mu )$ and (\ref{s1}) 
follows. 
\par\noindent
Again, for $\phi\in\mathcal P(K)$ and invoking Corollary 
\ref{new1} we have 
\[\left\|D_{A_0}^2\phi (x)\right\|^2_{HS}=\left\|\left(-A_0\right
)^{1/2}Q_{\infty}^{1/2}D^2\phi (x)Q_{\infty}^{1/2}\left(-A_0\right
)^{1/2}\right\|_{HS}^2\]
\begin{equation}=\left\|Q_{\infty}^{1/2}A_0Q_{\infty}^{1/2}D^2\phi 
(x)\right\|_{HS}^2=\left\|AQ_{\infty}D^2\phi (x)\right\|_{HS}^2=\frac 
14\left\|QD^2\phi (x)\right\|_{HS}^2,\label{c2}\end{equation}
and 
\[\left|A_0Q_{\infty}^{1/2}D\phi (x)\right|^2=\left\langle Q_{\infty}^{
1/2}A_0^2Q_{\infty}^{1/2}D\phi (x),D\phi (x)\right\rangle\]
\[=\left\langle AQ_{\infty}A^{*}D\phi (x),D\phi (x)\right\rangle 
=\frac 12\left\langle -QA^{*}D\phi (x),D\phi (x)\right\rangle\]
\begin{equation}=\frac 12\left|\left(-QA^{*}\right)^{1/2}D\phi (x
)\right|^2,\label{c3}\end{equation}
where the third equality follows from (\ref{65a}). Since the 
operators $D_{A_0}^2$ and $D_{A_0^2}$ are closable in $\mathcal W_
Q^{1,p}$ and 
$L^p(H,\mu )$ respectively, 
the above identities yield closability of of the operators 
$QD^2$ and $\left(-QA^{*}\right)^{1/2}D$ in $\mathcal W_Q^{1,p}$ and $
L^p(H,\mu )$ 
respectively, 
hence the corresponding 
Sobolev spaces 
$\mathcal W_{AQ}^{1,p}$ and $\mathcal W_Q^{2,p}$ are continuously embedded 
into $L^p(H,\mu )$. 
Finally, (\ref{c2}), (\ref{c3}) and Theorem \ref{td2} yield 
(\ref{s2}). 
\end{proof}
The Corollary below extends the result of \cite{lunardi} 
obtained by a completely different argument. 
\begin{corollary}
Assume that $\left(R_t\right)$ is a symmetric Ornstein-Uhlenbeck 
semigroup for a process in $\mathbb R^d$. Then for each $p\in (1,
\infty )$ 
\[\mbox{\rm dom}_p\left(\sqrt {-L}\right)=\mathcal W_I^{1,p}\quad\mbox{\rm and}
\quad\mbox{\rm dom}_p(L)=\mathcal W^{2,p}_I.\]
\end{corollary}
\begin{proof}
By Theorem \ref{t61a}(i) $Q$ has bounded inverse and the 
result follows immediately from Theorem \ref{t67}\ and 
the definition of Sobolev spaces. 
\end{proof}
The next corollary extends the results of \cite{shi}, \cite{dap}, 
\cite{regd} and \cite{dapgo}. 
\begin{corollary}
Assume that $A=A^{*}$ and $Q=(-A)^{-2\alpha}$, where $\alpha\ge 0$ and 
$\mbox{\rm tr}\left(-A\right)^{-1-2\alpha}<\infty$. Then $Q_{\infty}
=\frac 12\left(-A\right)^{-1-2\alpha}$ and for $p\in (1,\infty )$, 
\[\mbox{\rm dom}_p\left(\sqrt {-L}\right)=\mathcal W^{1,p}_{\left
(-A\right)^{-2\alpha}},\]
and
\[\mbox{\rm dom}_p(L)=\mathcal W_{\left(-A\right)^{-2\alpha}}^{2,
p}\cap\mathcal W_{(-A)^{1-2\alpha}}^{1,p}.\]
\end{corollary}
\begin{example}
In the previous Corollary assume that $H=L^2(0,1)$, and 
$A=\frac {\partial^2}{\partial\zeta^2}$ with zero Dirichlet boundary conditions and $
Q=I$. 
In this case $\mbox{\rm dom}\left(\sqrt {-A}\right)=H_0^1$ with equivalent norms. 
Then $\phi\in\mbox{\rm dom}_p(L)$ if and only if $D\phi (x)\in H_
0^1$ $\mu$-a.s., 
\[\int_H\left|D\phi (x)\right|^p_{H_0^1}\mu (dx)<\infty ,\]
$D^2\phi (x)$ is a Hilbert-Schmidt operator on $H$ $\mu$-a.s. and 
\[\int_H\left\|D^2\phi (x)\right\|_{HS}^p\mu (dx)<\infty .\]
\end{example}
\section{SPECTRAL GAP AND REGULARITY}
We start with the result which says that the spectral 
gap of the operator $L$ in $L^2(H,\mu )$ is the same as the 
spectral gap of $A_Q$ in $H$. 
\begin{theorem}\label{tg1}
Assume that the semigroup $\left(R_t\right)$ is symmetric. Then 
\begin{equation}\left|S_Q(t)h\right|\le e^{-\beta t}|h|,\quad h\in 
H,\label{er1}\end{equation}
if and only if  
\begin{equation}\left\|R_t\phi -\int_H\phi d\mu\right\|_2\le e^{-
\beta t}\left\|\phi\right\|_2,\quad\phi\in L^2(H,\mu ).\label{er2}\end{equation}
\end{theorem}
\begin{proof}
Assume that (\ref{er1}) holds and let $\Pi_0\phi =\int_H\phi d\mu$. 
By (\ref{gamma}) and the 
properties of the second quantization operator (Lemma 1c 
of \cite{fock}) we have 
\begin{equation}\left\|R_t-\Pi_0\right\|_2=\left\|S_Q(t)\right\|\le 
e^{-\beta t},\label{gamma1}\end{equation}
hence (\ref{er2}) is satisfied. The converse statement 
follows from (\ref{gamma1}) in the same way. 
\end{proof}
\begin{theorem}\label{t63} Assume that $\left(R_t\right)$ is symmetric. 
Then the following conditions are equivalent. 
\par\noindent
(i) $\left\|S_Q(t)\right\|=e^{-\beta t}$.
\par\noindent
(ii) $\left\|S_0(t)\right\|=e^{-\beta t}$.
\par\noindent
(iii) $\mbox{\rm im}\left(Q_t^{1/2}\right)=\mbox{\rm im}\left(Q_{
\infty}^{1/2}\right)$ for $t>0$. 
\par\noindent
(iv) $\mbox{\rm im}\left(Q_{\infty}^{1/2}\right)\subset\mbox{\rm im}\left
(Q^{1/2}\right)$.
\par\noindent
(v) The generator $L$ of $\left(R_t\right)$ satisfies the 
Logarithmic Sobolev Inequality: 
\[\int_H|\phi (x)|^2\log\,|\phi (x)|\,\mu (dx)\le\frac 2{\beta}\left
\langle -L\phi ,\phi\right\rangle +\|\phi\|^2\log\,\|\phi\|\]
\par\noindent
(vi) $\left(R_t\right)$ is hypercontractive from $L^p(H,\mu )$ to 
$L^q(H,\mu )$ for all $p,q$ such that  
\[\mbox{\rm $1<p<q\le 1+(p-1)e^{2\beta t}$.}\]
\end{theorem}\ 
\begin{proof}
$(i)\Leftrightarrow (ii)$. This follows immediately from (\ref{new3}). 
$(ii)\Leftrightarrow (iii)$ It is enough to 
recall that by Proposition 2b) in \cite{reg} (iii) is equivalent to $\left
\|S_0(t)\right\|<1$ 
for $t>0$. Then 
the symmetry of the 
semigroup $\left(S_0(t)\right)$ implies that for a certain $\beta 
>0$ we have 
$\left\|S_0(t)\right\|\le e^{-\beta t}$ for all $t\ge 0$. The same result in \cite{reg} 
shows that (ii) implies (iii). 
\par\noindent
$(ii)\Leftrightarrow (iv)$ By Proposition B.1 in \cite{dz1} (iv) holds if and 
only if there exists $a>0$ such that 
\begin{equation}\left|Vx\right|\ge a|x|,\quad x\in Q_{\infty}^{1/
2}(H).\label{new4}\end{equation}
By Lemma \ref{fock}\ (\ref{new4}) is equivalent to the 
condition
\[\left\langle A_0^{*}x,x\right\rangle\le -\frac 12a^2|x|^2,\quad 
x\in\mbox{\rm dom}\left(A_0^{*}\right),\]
and the last inequality is equivalent to (ii). 
\par\noindent
$(iii)\Leftrightarrow (v)\Leftrightarrow (vi)$ 
By Theorem 2 in \cite{fock} (iii) and (vi) are equivalent 
and by 
\cite{gross1}, see also \cite{rothaus}, (v) and (vi) are 
equivalent for symmetric semigroups. 
\end{proof}
\begin{remark}\label{tnew3}
Let us recall that for a finite Borel measure $\nu$ we denote 
by $\left\|\nu\right\|_{var}$ the variation norm of $\nu$ and the measure $
\nu R_t$ is 
defined by the formula 
\[\nu R_t(B)=\int_HR_tI_B(x)\nu (dx).\]
Assume that $\left(R_t\right)$ is symmetric. Then (\ref{er2}) (or any 
of the conditions of Theorem \ref{t63}) holds if and only if 
for each probability measure $\nu$ on $H$ such that $\nu\ll\mu$ 
and $\frac {d\nu}{d\mu}\in L^2(H,\mu )$ there exists $C_{\nu}<\infty$ such that 
\[\left\|\nu R_t-\mu\right\|_{var}\le C_{\nu}e^{-\beta t},\quad t
\ge 0.\]
This fact follows immediately from 
Theorem \ref{tg1} and the result in \cite{chen}, see 
also \cite{roberts}. 
\end{remark}
\begin{remark}
Let $\left(R_t\right)$ be a symmetric Ornstein-Uhlenbeck semigroup.
\par\noindent
(a) If $H=\mathbb R^d$ then by Theorem \ref{t61a}(i) $Q$ is 
boundedly invertible, hence (\ref{er1}) holds if and only if 
$\left(S(t)\right)$ is exponentially stable. Then by Theorem \ref{tg1}\ 
condition (\ref{er2}) is equivalent to the exponential 
stability of $\left(S(t)\right)$. In $H=\mathbb R^d$ the latter follows 
from Hypothesis \ref{H} (see Remark 2.6), hence 
 (\ref{er2}) always holds. 
\par\noindent
(b) If $\mbox{\rm dim}(H)=\infty$, then properties from (a) are not true 
in general. In Example 1 of Section 6 the 
Ornstein-Uhlenbeck semigroup is symmetric but it does 
not satisfy (\ref{er2}). In Example 2 the semigroup $\left(S(t)\right
)$ 
is not stable but (\ref{er2}) still holds. 
\end{remark}
The next corollary provides characterization of the 
Reproducing Kernel Hilbert Space $Q_{\infty}^{1/2}(H)$ of the invariant 
measure $\mu$. It follows immediately from Proposition 
\ref{trkhst} and Theorem   
\ref{t63}. Let us recall, that in the 
general nonsymmetric case the space $Q_{\infty}^{1/2}(H)$ is 
explicitly characterized as an interpolation space $D_A\left(\frac 
12,2\right)$ 
only if $\left(S(t)\right)$ is analytic and $Q$ is boundedly invertible, 
see Remark B.8 in \cite{dz1}
\begin{corollary}\label{trkhs}
Assume that $\left(R_t\right)$ is symmetric and $\left(S_Q(t)\right
)$ is 
exponentially stable. Then 
\[Q_{\infty}^{1/2}(H)=Q^{1/2}\left(\mbox{\rm dom}\sqrt {-A_Q}\right
).\]
\end{corollary}
\begin{theorem}\label{comp}
Assume that $\left(R_t\right)$ is symmetric. Then 
$\left(R_t\right)$ is a compact semigroup in $L^2(H,\mu )$ if and only if $\left
(S_Q(t)\right)$ is 
a compact and exponentially stable semigroup in $H$. 
\end{theorem}
\begin{proof}\ If $\left(R_t\right)$ is compact then by Proposition 2 in \cite{fock} and 
Theorem \ref{tnew1}\ 
for each $t>0$ $\left\|S_Q(t)\right\|<1$ and the semigroup $\left
(S_Q(t)\right)$ is 
compact. By the symmetry of $\left(S_Q(t)\right)$ the former implies that $\left
(S_Q(t)\right)$ is 
exponentially stable. If $\left(S_Q(t)\right)$ is compact and 
exponentially stable then $\left\|S_Q(t)\right\|<1$, hence $\left
(R_t\right)$ is compact by 
the result in \cite{fock}. 
\end{proof}
\begin{remark}\label{tangent} 
If $\left(R_t\right)$ is symmetric and compact then by Theorem 
\ref{t67}\ the embedding of 
$W_Q^{1,p}$ into $L^p(H,\mu )$ is compact for each $p\in (1,\infty 
)$. 
\end{remark}
The next result extends some results of Da Prato, see 
for example \cite{monotone}, where the estimate 
(\ref{grad}) is proved for the case $A=A^{*}$. In the 
theorem below we use the notation 
\[\|\phi\|_{\infty}=\mbox{\rm ess\,}\sup|\phi (x)|,\]
for $\phi\in L^{\infty}(H,\mu )$. 
\begin{theorem}\label{tgrad}
Assume that $\left(R_t\right)$ is symmetric and $\left(S_Q(t)\right
)$ is 
exponentially stable.  
Then $R_t$ is a bounded operator from $L^p(H,\mu )$ into 
$\mathcal W_Q^{1,p}$ for each 
$p\in (1,\infty )$ and $t>0$ and there exists $c(p)<\infty$ such that 
\begin{equation}\left\|Q^{1/2}DR_t\phi\right\|_p\le\frac {c(p)}{\sqrt {
t}}\left\|\phi\right\|_p.\label{grad}\end{equation}
Moreover, (\ref{grad})  still holds for $p=\infty$ for all 
bounded Borel functions $\phi$ and with the operator 
$\left(Q^{1/2}D,\mathcal W_Q^{1,p}\right)$ taken for an arbitrary $
p\in (1,\infty )$. 
\end{theorem}
\begin{proof}
For $p\in (1,\infty )$ the estimate (\ref{grad}) follows 
immediately from (\ref{s1}) and properties of analytic 
semigroups but we need another argument 
for the case $p=\infty$. Note that $V^{*}=Q_{\infty}^{-1/2}Q^{1/2}$ and 
$\mbox{\rm dom}\left(V^{*}\right)=\left\{x\in H:Q^{1/2}x\in Q_{\infty}^{
1/2}(H)\right\}$ is dense in $H$ since $V$ is closable. For $x\in\mbox{\rm dom}\left
(V^{*}\right)$
\[S(t)Q^{1/2}x=Q^{1/2}_{\infty}S_0(t)V^{*}x\in\mbox{\rm im}\left(
Q_t^{1/2}\right),\]
since $\mbox{\rm im}\left(Q_t^{1/2}\right)=\mbox{\rm im}\left(Q_{
\infty}^{1/2}\right)$ by Theorem \ref{t63}. Hence 
the operator 
\[Q_t^{-1/2}S(t)Q^{1/2}=Q_t^{-1/2}Q_{\infty}^{1/2}S_0(t)V^{*},\]
with the domain $\mbox{\rm dom}\left(V^{*}\right)$ is densely defined and since 
$V^{*}=\sqrt {-2A_0}U^{*}$ (see the proof of Theorem \ref{tnew1}) and $
Q_t^{-1/2}Q_{\infty}^{1/2}$ is bounded, it extends to a 
bounded operator on $H$. We will show that 
\begin{equation}\left\|Q_t^{-1/2}S(t)Q^{1/2}\right\|\le\frac 1{\sqrt {
t}}.\label{e1}\end{equation}
To this end note that for $h\in H$ 
\[Q_th=\int_0^tS(s)QS^{*}(s)hds=Q^{1/2}\int_0^tS_Q(2s)dsQ^{1/2}h,\]
hence for $x\in Q^{1/2}(H)$
\begin{equation}Q^{-1/2}Q_tQ^{-1/2}x=-\frac 12A_Q^{-1}\left(I-S_Q
(2t\right)x.\label{aaa}\end{equation}
It is easy to see from (\ref{67a}) that $\left(S_Q^{*}(t)\right)$ defines a $
C_0$-semigroup in 
$H_Q$ (see Remark \ref{t61b}) and since $\left(S_Q(t)\right)$ is 
symmetric we obtain from (\ref{aaa})  
for $x\in\mbox{\rm dom}\left(A_Q\left|H_Q\right.\right)$, the domain of the part of $
A_Q$ in $H_Q$, 
\[Q^{1/2}Q_t^{-1}Q^{1/2}x=-2A_Q\left(I-S_Q(2t)\right)^{-1}x.\]
Since $Q_t^{-1/2}S(t)Q^{1/2}=Q_t^{-1/2}Q^{1/2}S_Q(t)$ we obtain 
\[\left|Q_t^{-1/2}S(t)Q^{1/2}h\right|^2=\left|Q_t^{-1/2}Q^{1/2}S_
Q(t)h\right|^2\]
\[=2\left|\sqrt {-A_Q\left(I-S_Q(2t)\right)^{-1}S_Q(2t)}h\right|^
2.\]
By the Functional Calculus for selfadjoint operators 
\[\left\|A_Q\left(I-S_Q(2t)\right)^{-1}S_Q(2t)\right\|\le\sup_{\lambda 
>0}\left(\frac {\lambda e^{-2\lambda t}}{1-e^{-2\lambda t}}\right
)\le\frac 1{2t},\]
hence (\ref{e1}) holds. Fix $p\in (1,\infty )$ and let $D_I$ denote the 
closure in $L^p(H,\mu )$ of the Malliavin gradient (see below 
(\ref{s2})). By (ii) of Theorem \ref{t63}$ $ $A_0$ is boundedly 
invertible and by the argument in the proof of Theorem 
\ref{tnew1}\ so is  
$\bar {V}$, hence the operator $\bar {V}D_I$ with its maximal 
domain is closed in $L^p(H,\mu )$. Since $VD_I\phi =Q^{1/2}D\phi$ for 
$\phi\in\mathcal P(K)$ we conclude that $\bar {V}D_I\supset Q^{1/
2}D$. Let $\phi$ be 
bounded. By the first part of the theorem 
$R_t\phi\in\mathcal W_Q^{1,p}$, hence 
\begin{equation}Q^{1/2}DR_t\phi (x)=\bar {V}D_IR_t\phi (x),\quad\mu 
-a.s.\label{47a}\end{equation}
By Theorem 1 in \cite{reg}, condition (iii) of Theorem 
\ref{t63} implies that for a bounded Borel $\phi$ and $x\in H$, 
$D_IR_t\phi (x)$ exists as a Fr\'echet derivative in the direction 
$Q^{1/2}_{\infty}(H)$ and 
\begin{equation}\left\langle D_IR_t\phi (x),h\right\rangle =\int_
H\left\langle\Lambda^{*}(t)S_0(t)h,Q_t^{-1/2}y\right\rangle\phi (
S(t)x+y)\,\mu_t(dy),\label{reg}\end{equation}
where $\Lambda^{*}(t)S_0(t)=Q_t^{-1/2}S(t)Q_{\infty}^{1/2}$. 
Fix $x\in H$,  
such that (\ref{47a}) holds. Then for $h\in\mbox{\rm dom}\left(V^{
*}\right)$ 
(\ref{reg}) yields 
\[\left\langle Q^{1/2}DR_t\phi (x),h\right\rangle =\left\langle D_
IR_t\phi (x),V^{*}h\right\rangle\]
\[=\int_H\left\langle Q_t^{-1/2}S(t)Q^{1/2}h,Q_t^{-1/2}y\right\rangle
\phi\left(S(t)x+y\right)\mu_t(dy).\]
Therefore by (\ref{e1}) 
\[\left|\left\langle Q^{1/2}DR_t\phi (x),h\right\rangle\right|\le\sqrt {\frac 
2{\pi}}\left\|Q_t^{-1/2}S(t)Q^{1/2}\right\|\left\|\phi\right\|_{\infty}
|h|\]
\[\le\sqrt {\frac 2{\pi}}\frac 1{\sqrt t}\left\|\phi\right\|_{\infty}
|h|,\quad h\in\mbox{\rm dom}\left(V^{*}\right).\]
Since $\mbox{\rm dom}\left(V^{*}\right)$ is dense in $H$ we obtain (\ref{grad}) for 
$p=\infty$. 
\end{proof}
\begin{remark}\label{tangent1}
It follows from Remark \ref{t61b}\ and Theorem \ref{tg1}\ 
that (\ref{er2}) holds if and only 
\begin{equation}\left\|S_r(t)\right\|\le e^{-\beta t}.\label{rr2}\end{equation}
Moreover, $\left(R_t\right)$ is compact in $L^p(H,\mu )$, $p\in (
1,\infty )$ if and 
only if (\ref{rr2}) holds and $\left(S_r(t)\right)$ is compact in $
H_Q$. 
Finally, if (\ref{rr2}) holds then 
\[Q_{\infty}^{1/2}(H)=\mbox{\rm dom}\left(\sqrt {-A_r}\right),\]
by Remark \ref{t61b} and Theorem \ref{t63}.
\end{remark}
\section{HILBERT-SCHMIDT CASE}
We start with necessary and sufficient conditions for the 
semigroup $\left(R_t\right)$ to be Hilbert-Schmidt. This case was 
studied in \cite{reg}, where conditions for the 
Hilbert-Schmidt property were given in terms of the 
semigroup $\left(S_0(t)\right)$, hence rather difficult to apply in 
special cases. 
\begin{theorem}\label{ths1}
Assume that $\left(R_t\right)$ is symmetric. 
\par\noindent
(a) The following 
conditions are equivalent. 
\par\noindent\ 
(i) $\left(R_t\right)$ is a 
Hilbert-Schmidt semigroup in $L^2(H,\mu )$. 
\par\noindent
(ii) $\left(S_Q(t)\right)$ is a Hilbert-Schmidt and exponentially 
stable semigroup in $H$. 
\par\noindent
(b) Moreover,  $\mu\left(Q^{1/2}(H)\right)=1$ if and only if $\left
(S_Q(t)\right)$ is a 
Hilbert-Schmidt semigroup and 
\begin{equation}\int_0^{\infty}\left\|S_Q(t)\right\|_{HS}^2dt<\infty 
.\label{hj}\end{equation}
\end{theorem}
\begin{proof}
(a) In view of Proposition 2 in \cite{fock} the proof is completely analogous to the proof of 
Theorem \ref{comp} and therefore omitted. 
\par\noindent
(b) Assume that (\ref{hj}) holds. 
Let $H_Q$ denote the space defined in 
Remark \ref{t61b} and similarly, let $H_0=Q_{\infty}^{1/2}(H)$ be endowed 
with the norm $\left|x\right|_0=\left|Q_{\infty}^{-1/2}x\right|$. 
By Theorem \ref{t63} 
$Q_{\infty}^{1/2}(H)\subset Q^{1/2}(H)$, hence the corresponding imbedding 
$i:H_0\to H_Q$ is 
continuous. It is easy to check, that $\mu\left(Q^{1/2}(H)\right)
=1$ if and 
only if $i$ is a Hilbert-Schmidt operator, see also pp. 48-50 
of \cite{bogachev}. Let $\left\{e_k:k\ge 1\right\}$ be 
a CONS in $H$. Then $\left\{Q^{1/2}_{\infty}e_k:k\ge 1\right\}$ is a CONS in $
H_0$ and 
\[\sum_{k=1}^{\infty}\left|iQ^{1/2}_{\infty}e_k\right|^2_Q=\sum_{
k=1}^{\infty}\left|Q^{-1/2}Q_{\infty}^{1/2}e_k\right|^2.\]
Hence it is 
enough to show that the operator $Q^{-1/2}Q_{\infty}^{1/2}$ is 
Hilbert-Schmidt. 
For $x\in Q^{1/2}(H)$ we have 
\[Q_{\infty}Q^{-1/2}x=\int_0^{\infty}S(s)QS^{*}(s)Q^{-1/2}xds=Q^{
1/2}\int_0^{\infty}S_Q(2s)xds.\]
\begin{equation}=-\frac 12Q^{1/2}A_Q^{-1}x.\label{50}\end{equation}
Hence 
\[Q^{-1/2}Q_{\infty}Q^{-1/2}x=-\frac 12A_Q^{-1}x,\quad x\in Q^{1/
2}(H)\]
and thereby 
\begin{equation}Q^{-1/2}Q_{\infty}^{1/2}\overline {Q_{\infty}^{1/
2}Q^{-1/2}}=-\frac 12A_Q^{-1}.\label{piatek1}\end{equation}
Since $A_Q^{-1}$ is nuclear, (\ref{piatek1}) yields the 
Hilbert-Schmidt property of $Q^{-1/2}Q_{\infty}^{1/2}$ and thereby 
$\mu\left(Q^{1/2}(H)\right)=1$. 
\par\noindent
Conversely, assume that $\mu\left(Q^{1/2}(H)\right)=1$. Then 
it follows from the properties of 
Gaussian measures (see for example Theorem 2.5.8 in 
\cite{bogachevbook}) that 
$Q^{1/2}_{\infty}(H)\subset Q^{1/2}(H)$, hence by Theorem \ref{t63} $\left
(S_Q(t)\right)$ is 
exponentially stable. Consequently, $A_Q^{-1}$ is bounded and 
(\ref{50}) holds for $x\in Q^{1/2}(H)$, which implies 
(\ref{piatek1}). Since by (iii), $Q^{-1/2}Q_{\infty}^{1/2}$ is a 
Hilbert-Schmidt operator, it follows from (\ref{piatek1}) 
that $A_Q^{-1}$ is a nuclear operator. Hence, $S_Q(t)$ is Hilbert-Schmidt 
for each $t>0$ and 
\[\int_0^{\infty}\left\|S_Q(t)\right\|_{HS}^2ds=\frac 12\mbox{\rm tr}\left
(-A_Q^{-1}\right)<\infty .\]
\end{proof}
\begin{lemma}\label{tou1}
Let $Z$ be a solution to (\ref{01}) and assume that the 
corresponding Ornstein-Uhlenbeck semigroup is 
symmetric. Moreover, assume that (\ref{hj}) holds. 
Let $\tilde {Z}(\cdot ,x)$ denote a solution to the equation  
\begin{equation}\left\{\begin{array}{l}
d\tilde {Z}=A_Q\tilde {Z}dt+dW,\\
\tilde {Z}(0,x)=x\in H.\end{array}
\right.\label{ou1}\end{equation}
Then $Q^{-1/2}Z(t,x)=\tilde {Z}\left(t,Q^{-1/2}x\right)$ for each $
x\in Q^{1/2}(H)$. 
\end{lemma}
\begin{proof}\ 
By assumption the stochastic integral 
\[\int_0^tS_Q(t-s)dW(s)\]
is well defined. Moreover, 
since $Q^{-1/2}$ is closed in $H$ it is enough to note that for $
x\in Q^{1/2}(H)$  
\[Q^{-1/2}Z(t,x)=Q^{-1/2}S(t)x+Q^{-1/2}\int_0^tS(t-s)Q^{1/2}dW(s)\]
\[=S_Q(t)Q^{-1/2}x+\int_0^tS_Q(t-s)dW(s).\]
\end{proof}
Let $u(t,x)=E\phi\left(Z(t,x)\right)$ with $\phi\in L^p(H,\mu )$. We will show that 
for each $\phi\in L^p(H,\mu )$ the function $u$ satisfies almost 
everywhere an appropriate version of 
the following Backward Kolmogorov Equation 
\begin{equation}\left\{\begin{array}{l}
\frac {\partial u}{\partial t}(t,x)=\frac 12\mbox{\rm tr}\left(QD^
2u(t,x)\right)+\left\langle x,A^{*}Du(t,x)\right\rangle ,\\
u(0,x)=\phi (x).\end{array}
\right.\label{31}\end{equation}
For $\psi :H\to K$, where $K$ is a Banach space, let $D^Q\psi (x)$ 
denote the Fr\'echet derivative in the direction of the 
space  
$Q^{1/2}(H)$ of 
$\psi$ at the point $x\in H$. This means that $D^Q\psi (x)$ is a unique 
element of $\mathcal L(H,K)$ (note that $H=H^{*}$) such that 
\[\lim_{H\ni h\to 0}\frac {\psi\left(x+Q^{1/2}h\right)-\psi (x)-D^
Q\psi (x)h}{|h|}=0.\]
We write $\left(D^Q\right)^2\psi (x)=D^Q\left(D^Q\psi\right)(x)$. 
\begin{theorem}
Assume that $\left(R_t\right)$ is symmetric and the semigroup 
$\left(S_Q(t)\right)$ satisfies (\ref{hj}). Let $v(t,x)=u\left(t,
Q^{1/2}x\right)$. 
Then the following 
holds. 
\item{(a)} 
$v\in C^{1,2}\left((0,\infty )\times H,\mathbb R\right)$. 
\item{(b)} The functions  
$(t,x)\to\left\langle x,A_QD^Qu\left(t,Q^{1/2}x\right)\right\rangle$ and $
(t,x)\to\mbox{\rm tr}\left(\left(D^Q\right)^2u\left(t,Q^{1/2}x\right
)\right)$ are well defined and 
continuous on $(0,\infty )\times H$. 
\item{(c)} For every $t>0$ and 
$y\in Q^{1/2}(H)$ the function $u$ satisfies the following 
version of  
(\ref{31}): 
\begin{equation}\left\{\begin{array}{l}
\frac {\partial u}{\partial t}(t,y)=\frac 12\mbox{\rm tr}\left(\left
(D^Q\right)^2u(t,y)\right)+\left\langle Q^{-1/2}y,A_QD^Qu(t,y)\right
\rangle ,\\
u(0,y)=\phi (y).\end{array}
\right.\label{31a}\end{equation}
\end{theorem}
\begin{proof}
Let 
\[\tilde {R}_t\phi (x)=\mathbb E\phi\left(\tilde Z(t,x)\right).\]
Then we have for any $\phi\in L^p(H,\mu )$ 
\[v(t,x)=R_t\phi\left(Q^{1/2}x\right)=\mathbb E\phi\left(Q^{1/2}\tilde 
Z(t,x)\right)\]
\[=\mathbb E\tilde{\phi}\left(\tilde Z(t,x)\right)=\tilde {R}_t\tilde{
\phi }(x),\]
where $\tilde{\phi }(x)=\phi\left(Q^{1/2}x\right)$. Clearly, 
$\tilde{\phi}\in L^p(H,\mu )$. Moreover, by Corollary 9.22 in \cite{dz1} (or 
by (\ref{feller}) below) $\left(\tilde R_t\right)$ is strong Feller. Hence, by 
Theorem 5 in \cite{reg} 
the function $v(t,x)=\tilde {R}_t\tilde{\phi }(x)$ satisfies the 
following conditions. 
\par\noindent
(i) $v\in C^{1,2}\left((0,\infty )\times H,\mathbb R\right)$. 
\par\noindent
(ii) The functions $(t,x)\to A_QDv(t,x)$ and 
and $(t,x)\to\mbox{\rm tr}\left(D^2v(t,x)\right)$ are well defined and continuous 
on $(0,\infty )\times H$. 
\par\noindent
(iii) For every $t>0$ and $x\in H$ 
\begin{equation}\frac {\partial v}{\partial t}(t,x)=\frac 12\mbox{\rm tr}\left
(D^2v(t,x)\right)+\left\langle x,A_QDv(t,x)\right\rangle .\label{new4a}\end{equation}
Since $v$ is Fr\'echet differentiable for each 
$(t,x)\in (0,\infty )\times H$, the very definition of $v$ implies that 
\[\lim_{h\to 0}\frac {v(t,x+h)-v(t,x)-Dv(t,x)h}{|h|}\]
\[=\lim_{h\to 0}\frac {u\left(t,Q^{1/2}x+Q^{1/2}h\right)-u\left(t
,Q^{1/2}h\right)-Dv(t,x)h}{|h|}=0.\]
Hence, there exists $D^Qu\left(t,Q^{1/2}h\right)=Dv(t,x)$. Analogously, 
\[\left(D^Q\right)^2u\left(t,Q^{1/2}x\right)=D^2v(t,x),\quad (t,x
)\in (0,\infty )\times H.\]
Therefore, (b) follows from (ii) and (\ref{new4a}) yields 
\begin{equation}\frac {\partial u}{\partial t}\left(t,Q^{1/2}x\right
)=\frac 12\mbox{\rm tr}\left(\left(D^Q\right)^2\left(t,Q^{1/2}x\right
)\right)+\left\langle x,A_QD^Q\left(t,Q^{1/2}x\right)\right\rangle 
,\label{ddd}\end{equation}
for $t>0$ and $x\in H$. Putting $y=Q^{1/2}x$ in (\ref{ddd}) we 
obtain (\ref{31a}). 
\end{proof}
Let us recall that $\left(R_t\right)$ is strongly Feller if the 
function $R_t\phi$ is continuous for each $t>0$ and each bounded measurable 
function $\phi$. It was shown in \cite{sf} that the strong 
Feller property holds if and only if condition (\ref{sf}) 
below is satisfied, which is not easy to check in general. 
\begin{corollary}\label{tfeller}
Assume that $\left(R_t\right)$ is symmetric. Then $\left(R_t\right
)$ is strong Feller 
if and only if for each $t>0$ 
\begin{equation}\mbox{\rm im}\left(S(t)\right)\subset Q^{1/2}\left
(\mbox{\rm dom}\left(\sqrt {-A_Q}\right)\right).\label{feller}\end{equation}
If $\left(R_t\right)$ is strongly Feller then $\left(S_Q(t)\right
)$ is exponentially 
stable and of Hilbert-Schmidt type. 
\end{corollary}
\begin{proof}
By \cite{sf} (see also \cite{dz1}) $\left(R_t\right)$ is strongly Feller 
if and only if 
\begin{equation}\mbox{\rm im}\left(S(t)\right)\subset\mbox{\rm im}\left
(Q_t^{1/2}\right),\quad t>0.\label{sf}\end{equation}
By Proposition \ref{trkhst}\ 
\[\mbox{\rm im}\left(Q_t^{1/2}\right)=Q^{1/2}\left(\mbox{\rm dom}\left
(\sqrt {-A_Q}\right)\right),\quad t>0,\]
and consequently, (\ref{feller}) is equivalent to the strong 
Feller property of $\left(R_t\right)$. Let $\left(R_t\right)$ be strongly Feller. 
Then by \cite{dz1} (see also Proposition 3 in \cite{reg}) 
$\mbox{\rm im}\left(Q_t^{1/2}\right)=\mbox{\rm im}\left(Q_{\infty}^{
1/2}\right)$. Hence, by Theorem \ref{t63}\ $\left(S_Q(t)\right)$ 
is exponentially stable and by (\ref{sf}) $Q_{\infty}^{-1/2}S(t)$ is 
bounded for $t>0$. This implies that for $t>0$ the 
operator $S_0(t)=\left(Q_{\infty}^{-1/2}S(t)\right)Q_{\infty}^{1/
2}$ is of Hilbert-Schmidt 
type and by (\ref{68}) so is $S_Q(t)$ for all $t>0$. 
\end{proof}
\begin{remark}\label{trkhs1}
Remark \ref{trkhs}\ and Corollary \ref{tfeller}\ imply that 
the semigroup $\left(R_t\right)$ is strongly Feller if and only if 
\[\mbox{\rm im}\left(S(t)\right)\subset\mbox{\rm dom}\left(\sqrt {
-A_r}\right),\quad t>0.\]
\end{remark}
\section{EXAMPLES}
\subsection{Example 1}
The example below was introduced in \cite{impan} and 
later studied in a more general framework 
in \cite{jan}. Let $\left\{e_k:\,k\ge 1\right\}$ be a CONS in $H$ and let 
\[Qe_k=\frac 1{k^3}e_k\quad\mbox{\rm and}\quad Ae_k=-\frac 1ke_k.\]
Then $S(t)=e^{tA}$ and 
$\|S(t)\|=1$ for all $t\ge 0$. We have also 
\[Q_t=\frac 12A^2\left(I-e^{2tA}\right),\quad Q_{\infty}=\frac 12
A^2,\quad\mbox{\rm and}\quad A_Q=A.\]
We shall show that $\mbox{\rm im}\left(Q_t^{1/2}\right)$ is 
constant for all $t>0$ but $\mbox{\rm im}\left(Q_t^{1/2}\right)\neq\mbox{\rm $\mbox{\rm im}\left
(Q_{\infty}^{1/2}\right)$}$. Indeed, 
Proposition \ref{trkhst}\ yields
\[\mbox{\rm im}\left(Q_t^{1/2}\right)=Q^{1/2}(H)=A^{3/2}(H),\quad 
t>0,\]
while $Q_{\infty}^{1/2}(H)=A(H)$. 
It follows from Theorem \ref{t63} 
that $R_t$ is not hypercontractive for any $t>0$ and the 
generator $L$ of $\left(R_t\right)$ has no spectral gap. 
\par\ 
Using Theorem 
1a 
from \cite{reg} we find that there exists a bounded Borel function 
$\phi$, $x\in H$ and $h\in Q_{\infty}^{1/2}(H)$ such that the function $
t\to R_t\phi (x+th)$ 
is not continuous. 
\par
Let us recall that 
for noninteger $\alpha$ the space $\mathcal W^{\alpha ,2}_{Q_{\infty}}
(H)$ is defined by 
interpolation (for details, see 
\cite{fock}). It 
follows from Theorem 
4c in \cite{fock} that $R_t\left(L^2(H,\mu )\right)$ is not contained 
in $\mathcal W^{\alpha ,2}_{Q_{\infty}}(H)$ for 
any $\alpha >0$. Hence 
$\mbox{\rm dom}\,(L)$ is not contained in $W^{\alpha ,2}_{Q_{\infty}}
(H)$ for any $\alpha >0$. \hfill This 
fact can be also directly deduced from Theorem \ref{t67}\ 
which yields 
\[\mbox{\rm $\mbox{\rm dom}_2(L)=\mathcal W_{-A^3}^{2,2}\cap\mathcal 
W_{-A^4}^{1,2}$}.\]
\subsection{Example 2}
The stochastic heat equation (\ref{ouk}) in a weighted 
space $L^2\left(\mathbb R,\rho (\zeta )d\zeta\right)$ was considered in \cite{dz2} as 
an example of the Ornstein-Uhlenbeck process in a chaotic 
environment. Here we investigate some properties of the 
transition semigroup associated to (\ref{ouk}), using 
different methods. In particular, we improve some 
results from \cite{dz2}. 
\par
Let $H^{\kappa}=L^2\left(\mathbb R,\rho_{\kappa}(\zeta )d\zeta\right
)$, where $\rho_{\kappa}(\zeta )=e^{-\kappa |\zeta |}$ with 
$\kappa\ge 0$. In particular $H^0=L^2\left(\mathbb R\right)$. The scalar 
product and the norm in $H^{\kappa}$ will be denoted by $\left\langle
\cdot ,\cdot\right\rangle_{\kappa}$ 
and $|\cdot |_{\kappa}$ respectively. 
Fix $m>0$ and let $A^{(0)}=\Delta -mI$, where $\Delta$ is the Laplacian in 
$H^0$ and let $S^{(0)}(t)$ denote the semigroup on $H^0$ generated 
by $A^{(0)}$. Then $A^{(0)}$ is selfadjoint in $L^2(\mathbb R)$ and 
$\mbox{\rm dom}\sqrt {-A^{(0)}}=H^{1,2}(\mathbb R)$. The semigroup $\left
(S^{(0)}(t)\right)$ 
generated by $A^{(0)}$ has the property 
\begin{equation}\left\|S^{(0)}(t)\right\|=e^{-mt}.\label{f0}\end{equation}
Let $H_{\kappa}^{1,2}$ denote the space of functions $x\in H_{\kappa}$ such that 
the distributional derivative $x'\in H_{\kappa}$ and 
\[\left|x\right|^2_{\kappa ,1}=\left|x'\right|^2_{\kappa}+m|x|^2_{
\kappa}<\infty .\]
\begin{lemma}\label{tcoerc}
For any $\kappa ,m\ge 0$ there exist $\alpha >0$ and $\omega\in\mathbb 
R$ such 
that 
\begin{equation}\left\langle -A^{(0)}x,x\right\rangle_{\kappa}\ge
\alpha |x|^2_{\kappa ,1}+\omega |x|^2_k,\quad x\in C_0^{\infty}(\mathbb 
R).\label{coerc}\end{equation}
In particular, $A^{(0)}$ extends to a generator $A^{(\kappa )}$ of an analytic 
semigroup $\left(S^{(\kappa )}(t)\right)$ in $H^{\kappa}$ and 
\begin{equation}\left\|S^{(\kappa )}(t)\right\|\le e^{-\omega t}.\label{f1}\end{equation}
Moreover, if $\frac 14\kappa^2<m$ then $\omega >0$. 
\end{lemma}
\begin{proof}
For $x\in C_0^{\infty}(\mathbb R)$ we have 
\[\left\langle -A^{(0)}x,x\right\rangle_{\kappa}=\left\langle -A^{
(0)}x,\rho_{\kappa}x\right\rangle_0=\left\langle x',\left(\rho_kx\right
)^{'}\right\rangle_0+m|x|^2_{\kappa}\]
\[=\left|x'\right|_{\kappa}^2+\left\langle xx',\rho_k'\right\rangle_
0+m|x|^2_{\kappa}.\]
Since $\left|\rho_k'(\zeta )\right|=k\rho_k(\zeta )$ for $\zeta\neq 
0$ we obtain for any $\epsilon >0$: 
\[\left\langle -A^{(0)}x,x\right\rangle_{\kappa}\ge\left|x'\right
|_{\kappa}^2+m|x|^2_{\kappa}-\frac {k\epsilon}2\left|x'\right|_{\kappa}^
2-\frac k{2\epsilon}|x|^2_{\kappa}\]
\[=\left(1-\frac {k\epsilon}2\right)\left(\left|x'\right|^2_{\kappa}
+m|x|^2_{\kappa}\right)+\left(\frac {mk\epsilon}2-\frac k{2\epsilon}\right
)|x|^2_{\kappa}.\]
Hence, (\ref{coerc}) follows provided 
$k\epsilon <2$. The remaining part of the lemma follows easily 
from the Theorem of Lions (see p. 389 of \cite{dz1})
\end{proof}
We will consider equation (\ref{01}) written in a slightly 
different form 
\begin{equation}dZ=A^{(\kappa )}Zdt+JdW,\label{ouk}\end{equation}
where $W$ is standard cylindrical Wiener process on $H^{(0)}$ 
and $J:H^{(0)}\to H^{(\kappa )}$ is an embedding: $Jx=x$. Then $Q
=JJ^{*}$ 
and it is easy to check that 
\begin{equation}J^{*}x=\rho_{\kappa}x=Qx,\quad Q^{1/2}x=\rho_{\kappa}^{
1/2}x.\label{f2}\end{equation}
It was proved in \cite{dz2} that for any $\kappa >0$ and $m>0$ 
the solution (\ref{ouk}) is well defined in $H^{\kappa}$ and it 
admits an invariant measure $\mu =N\left(0,Q_{\infty}\right)$. Let $\left
(R_t\right)$ be the 
Ornstein-Uhlenbeck semigroup corresponding to (\ref{ouk}). 
\begin{proposition}\label{ex31}
For any $\kappa >0$ and $m>0$ the following holds. 
\par\noindent
(i) $\mbox{\rm ker}\left(Q_{\infty}\right)=\{0\}$ in $H^{\kappa}$. 
\par\noindent
(ii) $R_t=R_t^{*}$ in $L^2\left(H^{\kappa},\mu\right)$. 
\par\noindent
(iii) The semigroup $\left(R_t\right)$ 
satisfies all the statements of Theorem \ref{t63} 
with $\beta =m$. 
\end{proposition}
\begin{proof}
(i) Note that if (\ref{ouk}) has an invariant measure 
$\mu =N\left(0,Q_{\infty}\right)$ then 
\begin{equation}\mbox{\rm ker}\left(Q_{\infty}\right)\subset\mbox{\rm ker}\left
(Q\right).\label{f3}\end{equation}
Indeed, if $x\in\mbox{\rm ker}\left(Q_{\infty}\right)$ then 
\[0=\left\langle Q_{\infty}x,x\right\rangle =\int_0^{\infty}\left
|Q^{1/2}S^{*}(t)x\right|^2dt.\]
Hence, for a.a. $t\ge 0$, $Q^{1/2}S^{*}(t)x=0$ and by continuity 
$Q^{1/2}x=0$. Thus (i) follows from (\ref{f1}) and (\ref{f3}). 
\par\noindent
(ii) For $x\in H^0$ and $y\in H^{\kappa}$
\[\left\langle S^{(\kappa )}(t)x,y\right\rangle_{\kappa}=\left\langle 
S^{(0)}(t)x,\rho_{\kappa}y\right\rangle_0\]
\[=\left\langle x,S^{(0)}(t)\rho_{\kappa}y\right\rangle_0=\left\langle 
x,\rho_{\kappa}^{-1}S^{(0)}(t)\left(\rho_{\kappa}y\right)\right\rangle_{
\kappa},\]
and thereby 
\[\left(S^{(\kappa )}(t)\right)^{*}y=\rho_{\kappa}^{-1}S^{(0)}(t)\left
(\rho_{\kappa}y\right),\]
and by (\ref{f2}) 
\[Q\left(S^{(\kappa )}(t)\right)^{*}y=S^{(0)}(t)\left(\rho_{\kappa}
y\right)=S^{(\kappa )}(t)Qy.\]
Therefore, (ii) holds by Theorem \ref{t61}. 
\par\noindent
(iii) By Theorem \ref{t61a}\ and (\ref{f2}) we find that for 
$x\in H^{\kappa}$ 
\[\left|S^{(\kappa )}_Q(t)x\right|_{\kappa}^2=\int_{\mathbb R}\left
(\rho_{\kappa}^{-1/2}S^{(\kappa )}(t)\left(\rho_{\kappa}x\right)(
\zeta )\right)^2\rho_{\kappa}(\zeta )d\zeta\]
\[=\left|S^{(0)}(t)\left(\rho_{\kappa}^{1/2}x\right)\right|_0^2.\]
Then by (\ref{f0}) 
\begin{equation}\left\|S_Q^{(\kappa )}(t)\right\|_{\kappa}=e^{-mt}
,\label{f4}\end{equation}
and (iii) follows from Theorem (\ref{t63}).
\end{proof}
\begin{corollary}\label{ex32}
For any $\kappa >0$ and $m>0$ the following holds. 
\par\noindent
(i) For $\phi\in L^2\left(H^{\kappa},\mu\right)$ 
\[\left\|R_t\phi -\int_{H^{\kappa}}\phi d\mu\right\|_2=e^{-mt}\left
\|\phi\right\|_2.\]
(ii) $\left(R_t\right)$ is not strong Feller on $H^{\kappa}$. 
\par\noindent
(iii) 
\[Q_{\infty}^{1/2}\left(H^{\kappa}\right)=\mbox{\rm dom}\left(\sqrt {
-A^{(0)}}\right)=H^{1,2}\left(\mathbb R\right)=Q_t^{1/2}\left(H^{
\kappa}\right).\]
\end{corollary}
\begin{proof}
Part (i) follows from (\ref{f4}) and Theorem \ref{tg1}. 
Note that $H_Q=H^0$ and $\left(S^{(\kappa )}(t)\right)$ restricted to $
H_Q$ is 
isometrically isomorphic to $\left(S^{(0)}(t)\right)$. Since $\left
(S^{(0)}(t)\right)$ is 
not compact, the semigroup $\left(R_t\right)$ is not compact by 
Theorem \ref{comp}\ and Remark \ref{tangent}. Hence (ii) 
follows. Similarly, we obtain (iii) from Remark 
\ref{trkhs1}. 
\end{proof}
If $\frac 14\kappa^2<m$ then by (\ref{f1}) the semigroup $\left(S^{
(\kappa )}(t)\right)$ is 
exponentially stable, hence $\mu$ is a unique invariant 
measure for (\ref{ouk}). It was shown in \cite{dz2} that 
for $0<m<\frac 14\kappa^2$ there are infinitely many invariant 
measures for (\ref{ouk}). Below we improve this result. 
\begin{proposition}\label{ex33}
If $m\ge\frac 14\kappa^2$ then there exists a unique invariant measure. 
For $0<m<\frac 14\kappa^2$ there exists a family $\left\{\mu_{\lambda}
:\lambda\in\mathbb R\right\}$ 
of Borel probability measures on $H^{\kappa}$ such that for any 
$\lambda\in\mathbb R$ the following holds. 
\par\noindent
(i) The measure $\mu_{\lambda}$ is symmetrizing for $\left(R_t\right
)$.
\par\noindent
(ii) The Logarithmic 
Sobolev Inequality holds in $L^2\left(H^{\kappa},\mu_{\lambda}\right
)$:
\[\int_{H^{\kappa}}\phi^2(x)\log\left|\phi (x)\right|^2\mu (dx)\le\frac 
2m\left\langle -L\phi ,\phi\right\rangle_{\kappa}+\left\|\phi\right
\|^2_{\kappa}\log\left\|\phi\right\|_{\kappa}^2.\]
\par\noindent
(iii) If $\lambda_1\neq\lambda_2$ then $\mu_{\lambda_1}$ and $\mu_{
\lambda_2}$ are singular and in 
particular, $\mu_{\lambda}\perp\mu =\mu_0$ if $\lambda\neq 0$. 
\end{proposition}
\begin{proof}
(i) If $m>\frac 14\kappa^2$ then the uniqueness of the invariant measure 
follows from (\ref{f1}). For $m=\frac 14\kappa^2$, assume that for a 
certain $x\in\mbox{\rm dom}(A)$ the 
equation $A^{(\kappa )}x=0$ has a solution. Then $\Delta x=mx$, and 
therefore $x\in\mbox{\rm dom}\left((-\Delta )^n\right)$ for all $
n\ge 1$. Hence 
$x\in C^2(\mathbb R)$ and $x^{\prime\prime}=mx$, which is impossible for 
$x\in H^{\kappa}$. 
\par\noindent
(iii) Following \cite{dz2}, let $g(\zeta )=e^{\sqrt {m}\zeta}$. Then 
\begin{equation}A^{(\kappa )}g=0.\label{f5}\end{equation}
Therefore, $\mu_{\lambda}=N\left(\lambda g,Q_{\infty}\right)$, $\lambda
\in\mathbb R$, is also 
invariant for (\ref{ouk}), see \cite{dz2} for details. By 
(iii) of 
Corollary \ref{ex32}, $Q_{\infty}^{1/2}\left(H^{\kappa}\right)\subset 
H^0$. Since for $\lambda\neq 0$, 
$\lambda g\notin H^0$ we obtain (iii) from the Feldman-Hayek Theorem. 
\par\noindent
(i) Since $\mu_{\lambda}$ is invariant for (\ref{ouk}), $\left(R_
t\right)$ is a 
$C_0$-semigroup of contractions on $L^p\left(H^{\kappa},\mu_{\lambda}\right
)$ for 
$p\in [1,\infty )$. For $a\in H^{\kappa}$ let $T_a$ denote the shift operator on 
$H^{\kappa}$: $T_ax=a+x$. Note that for for a bounded Borel $\phi
\ge 0$ 
and $\mu_a=N\left(a,Q_{\infty}\right)$ we have 
\begin{equation}\int_{H^{\kappa}}\phi\left(T_ax\right)\mu (dx)=\int_{
H^{\kappa}}\phi (x)\mu_a(dx),\label{f6}\end{equation}
hence the map 
\begin{equation}\phi\to\phi\circ T_a:L^p\left(H^{\kappa},\mu_a\right
)\to L^p\left(H^{\kappa},\mu\right),\label{f6a}\end{equation}
is an isometric isomorphism. For $a=\lambda g$ (\ref{f5}) yields 
\[Z\left(t,a+x\right)=S^{(\kappa )}(t)a+Z(t,x)=a+Z(t,x),\]
which implies that 
\begin{equation}R_t\phi\left(T_ax\right)=R_t\left(\phi\circ T_a\right
)(x),\quad\phi\in L^1\left(H^{\kappa},\mu_{\lambda}\right).\label{f7}\end{equation}
Taking into account that $\left(R_t\right)$ is symmetric in $L^2\left
(H^{\kappa},\mu\right)$, 
by (\ref{f6}) and (\ref{f7}) we have for $\phi ,\psi\in L^2\left(
H^{\kappa},\mu_{\lambda}\right)$
\[\int_{H^{\kappa}}\psi (x)R_t\phi (x)\mu_{\lambda}(dx)=\int_{H^{
\kappa}}\phi\circ T_a(x)R_t\left(\psi\circ T_a\right)(x)\mu (dx)\]
\[=\int_{H^{\kappa}}\phi (x)R_t\psi (x)\mu_{\lambda}(dx),\]
which proves (i). 
\par\noindent
(ii) It follows from (iii) of Proposition 
\ref{ex31}, (\ref{f6a}) and (\ref{f7}) that the semigroup $\left(
R_t\right)$ 
is hypercontractive in each $L^p\left(H^{\kappa},\mu_{\lambda}\right
)$ and therefore (ii) 
follows from (i). 
\end{proof}
\begin{remark}\label{ex4}
If $m<\frac 14\kappa^2$ then the asymptotic behaviour of 
the semigroups $\left(S_Q^{(\kappa )}(t)\right)$ and $\left(S^{(\kappa 
)}(t)\right)$ on $H^{\kappa}$ is 
different: while the former is exponentially stable, the 
latter is not.  Actually, for $m<\frac 14\kappa^2$ 
\[\left\|S_Q^{(\kappa )}(t)\right\|\le e^{-mt},\quad\lim_{t\to\infty}\left
\|S^{(\kappa )}(t)\right\|\to\infty .\]
Indeed, let $x_{\alpha}(\zeta )=e^{\alpha\zeta}$. Then $x_{\alpha}
\in H^{\kappa}$ (but not 
to $H^0$), provided $|\alpha |<\frac 12\kappa$ and $\left(\Delta 
-mI\right)x_{\alpha}=\left(\alpha^2-m\right)x_{\alpha}$, hence 
$S^{(\kappa )}(t)x_{\alpha}=e^{t\left(\alpha^2-m\right)}x_{\alpha}$. Thereby, for $
m<\alpha^2<\frac 14\kappa^2$ we obtain 
$\left|S^{(\kappa )}(t)x_{\alpha}\right|_{\kappa}=e^{t\left(\alpha^
2-m\right)}\left|x_{\alpha}\right|_{\kappa}\to\infty$ for $t\to\infty$. 
\end{remark}
\begin{remark}
We show in part (iii) of Proposition \ref{ex33}\ that the 
invariant measures $\mu_{\lambda}$ must be mutually singular. Since 
the generator $L$ of $\left(R_t\right)$ in $L^2\left(H^{\kappa},\mu_{
\lambda}\right)$ can be associated 
to an irreducible symmetric Dirichlet form, this 
fact can be also deduced from  
\cite{roeckner}. For similar results for processes which 
are not associated to Dirichlet forms 
see \cite{ania}. 
\end{remark}
\par\medskip\noindent
\centerline{ACKNOWLEDGEMENT}
\par\noindent
The authors are indebted to the referee for the numerous 
suggestions which improved presentation of the paper. 

\end{document}